\documentclass[a4paper,11pt,twoside]{article}
\usepackage{amsmath,amssymb,amsthm}

\renewcommand{\k}{\Bbbk}
\newcommand{\N}{\mathbb{N}}
\newcommand{\Z}{\mathbb{Z}}
\newcommand{\R}{\mathbb{R}}

\DeclareMathOperator{\codim}{codim}

\DeclareMathOperator{\Lin}{Lin}
\DeclareMathOperator{\cop}{\triangle}
\newcommand{\tensor}{\otimes}
\DeclareMathOperator{\cou}{\epsilon}
\DeclareMathOperator{\antip}{S}
\newcommand{\act}{\triangleright}

\newcommand{\xprod}{\!\bowtie\!}
\newcommand{\lalg}[1]{\mathfrak{#1}}

\newcommand{\no}[1]{{:{#1}:}}
\DeclareMathOperator{\id}{id}
\DeclareMathOperator{\adl}{Ad_L}
\newcommand{\binomq}[2]{\genfrac{[}{]}{0pt}{}{#1}{#2}_q}
\newcommand{\cbp}{\ensuremath{C(B_+)}}
\newcommand{\cqbp}{\ensuremath{C_q(B_+)}}
\newcommand{\cbnp}{\ensuremath{C(B_{n+})}}
\newcommand{\ubp}{\ensuremath{U(\lalg{b_+})}}
\newcommand{\uqbp}{\ensuremath{U_q(\lalg{b_+})}}
\newcommand{\ubnp}{\ensuremath{U(\lalg{b}_{n+})}}
\DeclareMathOperator{\diff}{d}
\DeclareMathOperator{\fdiff}{\nabla}

\theoremstyle{plain}
\newtheorem{prop}{Proposition}[section]
\newtheorem{lem}[prop]{Lemma}
\newtheorem{thm}[prop]{Theorem}
\newtheorem{cor}[prop]{Corollary}

\title{Classification of Differential Calculi on \uqbp, Classical
Limits, and Duality}
\author{Robert Oeckl\\
Department of Applied Mathematics and Theoretical Physics,\\
University of Cambridge, Cambridge CB3 9EW, UK}
\date{DAMTP-1998-86\\
July 18,1998}
\begin{document}
\maketitle

\begin{abstract}
We give a complete classification of bicovariant first order
differential
calculi on the quantum enveloping algebra \uqbp{} which we view as the
quantum function algebra \cqbp{}. Here, $\lalg{b}_+$ is the Borel subalgebra
of $\lalg{sl}_2$. We do the same in the classical
limit $q\to 1$ and obtain a one-to-one correspondence in the finite
dimensional case. It turns out that the classification is essentially
given by finite subsets of the positive integers. We proceed to investigate
the classical limit from the dual point of view, i.e. with ``function
algebra'' \ubp{} and ``enveloping algebra'' \cbp{}. In this case there
are many more differential calculi than coming from the $q$-deformed
setting. As an application, we give the natural intrinsic
4-dimensional calculus of
$\kappa$-Minkowski space and the
associated formal integral.
\end{abstract}

\section{Introduction}

One of the fundamental ingredients in the theory of non-commutative or
quantum geometry is the notion of a differential calculus.
In the framework of quantum groups the natural notion
is that of a
bicovariant differential calculus as introduced by Woronowicz
\cite{Wor_calculi}. Due to the allowance of non-commutativity
the uniqueness of a canonical calculus is lost.
It is therefore desirable to classify the possible choices.
The most important piece is the space of one-forms or ``first
order differential calculus'' to which we will restrict our attention
in the following. (From this point on we will use the term
``differential calculus'' to denote a
bicovariant first order differential calculus).

Much attention has been devoted to the investigation of differential
calculi on quantum groups $C_q(G)$ of function algebra type for
$G$ a simple Lie group.
Natural differential calculi on matrix quantum groups were obtained by
Jurco \cite{Jur} and
Carow-Watamura et al.\
\cite{CaScWaWe}. A partial classification of calculi of the same
dimension as the natural ones
was obtained by
Schm\"udgen and Sch\"uler \cite{ScSc2}.
More recently, a classification theorem for factorisable
cosemisimple quantum groups was obtained by Majid \cite{Majid_calculi},
covering the general $C_q(G)$ case. A similar result was
obtained later by Baumann and Schmitt \cite{BaSc}.
Also, Heckenberger and Schm\"udgen \cite{HeSc} gave a
complete classification on $C_q(SL(N))$ and $C_q(Sp(N))$.

In contrast, for $G$ not simple or semisimple the differential calculi
on $C_q(G)$
are largely unknown. A particularly basic case is the Lie group $B_+$
associated with the Lie algebra $\lalg{b_+}$ generated by two elements
$X,H$ with the relation $[H,X]=X$. The quantum enveloping algebra
\uqbp{}
is self-dual, i.e.\ is non-degenerately paired with itself \cite{Drinfeld}.
This has an interesting consequence: \uqbp{} may be identified with (a
certain algebraic model of) \cqbp. The differential calculi on this
quantum group and on its ``classical limits'' \cbp{} and \ubp{}
will be the main concern of this paper. We pay hereby equal attention
to the dual notion of ``quantum tangent space''.

In section \ref{sec:q} we obtain the complete classification of differential
calculi on \cqbp{}. It turns out that (finite
dimensional) differential
calculi are characterised by finite subsets $I\subset\N$.
These
sets determine the decomposition into coirreducible (i.e.\ not
admitting quotients) differential calculi
characterised by single integers. For the coirreducible calculi the
explicit formulas for the commutation relations and braided
derivations are given.

In section \ref{sec:class} we give the complete classification for the
classical function algebra \cbp{}. It is essentially the same as in the
$q$-deformed setting and we stress this by giving an almost
one-to-one correspondence of differential calculi to those obtained in
the previous section. In contrast, however, the decomposition and
coirreducibility properties do not hold at all. (One may even say that
they are maximally violated). We give the explicit formulas for those
calculi corresponding to coirreducible ones.

More interesting perhaps is the ``dual'' classical limit. I.e.\ we
view \ubp{} as a quantum function algebra with quantum enveloping
algebra \cbp{}. This is investigated in section \ref{sec:dual}. It
turns out that in this setting we have considerably more freedom in
choosing a
differential calculus since the bicovariance condition becomes much
weaker. This shows that this dual classical limit is in a sense
``unnatural'' as compared to the ordinary classical limit of section
\ref{sec:class}. 
However, we can still establish a correspondence of certain
differential calculi to those of section \ref{sec:q}. The
decomposition properties are conserved while the coirreducibility
properties are not.
We give the
formulas for the calculi corresponding to coirreducible ones.

Another interesting aspect of viewing \ubp{} as a quantum function
algebra is the connection to quantum deformed models of space-time and
its symmetries. In particular, the $\kappa$-deformed Minkowski space
coming from the $\kappa$-deformed Poincar\'e algebra
\cite{LuNoRu}\cite{MaRu} is just a simple generalisation of \ubp.
We use this in section \ref{sec:kappa} to give
a natural $4$-dimensional differential calculus. Then we show (in a
formal context) that integration is given by
the usual Lesbegue integral on $\R^n$ after normal ordering.
This is obtained in an intrinsic context different from the standard
$\kappa$-Poincar\'e approach.

A further important motivation for the investigation of differential
calculi on
\ubp{} and \cbp{} is the relation of those objects to the Planck-scale
Hopf algebra \cite{Majid_Planck}\cite{Majid_book}. This shall be
developed elsewhere.

In the remaining parts of this introduction we will specify our
conventions and provide preliminaries on the quantum group \uqbp, its
deformations, and differential calculi.

\subsection{Conventions}

Throughout, $\k$ denotes a field of characteristic 0 and
$\k(q)$ denotes the field of rational
functions in one parameter $q$ over $\k$.
$\k(q)$ is our ground field in
the $q$-deformed setting, while $\k$ is the
ground field in the ``classical'' settings.
Within section \ref{sec:q} one could equally well view $\k$ as the ground
field with $q\in\k^*$ not a root of unity. This point of view is
problematic, however, when obtaining ``classical limits'' as
in sections \ref{sec:class} and \ref{sec:dual}.

The positive integers are denoted by $\N$ while the non-negative
integers are denoted by $\N_0$.
We define $q$-integers, $q$-factorials and
$q$-binomials as follows:
\begin{gather*}
[n]_q=\sum_{i=0}^{n-1} q^i\qquad
[n]_q!=[1]_q [2]_q\cdots [n]_q\qquad
\binomq{n}{m}=\frac{[n]_q!}{[m]_q! [n-m]_q!}
\end{gather*}
For a function of several variables (among
them $x$) over $\k$ we define
\begin{gather*}
(T_{a,x} f)(x) = f(x+a)\\
(\fdiff_{a,x} f)(x) = \frac{f(x+a)-f(x)}{a}
\end{gather*}
with $a\in\k$ and similarly over $\k(q)$
\begin{gather*}
(Q_{m,x} f)(x) = f(q^m x)\\
(\partial_{q,x} f)(x) = \frac{f(x)-f(qx)}{x(1-q)}\\
\end{gather*}
with  $m\in\Z$.

We frequently use the notion of a polynomial in an extended
sense. Namely, if we have an algebra with an element $g$ and its
inverse $g^{-1}$ (as
in \uqbp{}) we will mean by a polynomial in $g,g^{-1}$ a finite power
series in $g$ with exponents in $\Z$. The length of such a polynomial
is the difference between highest and lowest degree.

If $H$ is a Hopf algebra, then $H^{op}$ will denote the Hopf algebra
with the opposite product.

\subsection{\uqbp{} and its Classical Limits}
\label{sec:intro_limits}

We recall that,
in the framework of quantum groups, the duality between enveloping algebra
$U(\lalg{g})$ of the Lie algebra and algebra of functions $C(G)$ on the Lie
group carries over to $q$-deformations.
In the case of
$\lalg{b_+}$, the
$q$-deformed enveloping algebra \uqbp{} defined over $\k(q)$ as
\begin{gather*}
U_q(\lalg{b_+})=\k(q)\langle X,g,g^{-1}\rangle \qquad
\text{with relations} \\
g g^{-1}=1 \qquad Xg=qgX \\
\cop X=X\tensor 1 + g\tensor X \qquad
\cop g=g\tensor g \\
\cou (X)=0 \qquad \cou (g)=1 \qquad
\antip X=-g^{-1}X \qquad \antip g=g^{-1}
\end{gather*}
is self-dual. Consequently, it
may alternatively be viewed as the quantum algebra \cqbp{} of
functions on the Lie group $B_+$ associated with $\lalg{b_+}$.
It has two classical limits, the enveloping algebra \ubp{}
and the function algebra $C(B_+)$.
The transition to the classical enveloping algebra is achieved by
replacing $q$
by $e^{-t}$ and $g$ by $e^{tH}$ in a formal power series setting in
$t$, introducing a new generator $H$. Now, all expressions are written in
the form $\sum_j a_j t^j$ and only the lowest order in $t$ is kept.
The transition to the classical function algebra on the other hand is
achieved by setting $q=1$.
This may be depicted as follows:
\[\begin{array}{c @{} c @{} c @{} c}
& \uqbp \cong \cqbp && \\
& \diagup \hspace{\stretch{1}} \diagdown && \\
 \begin{array}{l} q=e^{-t} \\ g=e^{tH} \end{array} \Big| _{t\to 0} 
 && q=1 &\\
 \swarrow &&& \searrow \\
 \ubp & <\cdots\textrm{dual}\cdots> && \cbp
\end{array}\]
The self-duality of \uqbp{} is expressed as a pairing
$\uqbp\times\uqbp\to\k$
with
itself:
\[\langle X^n g^m, X^r g^s\rangle =
 \delta_{n,r} [n]_q!\, q^{-n(n-1)/2} q^{-ms}
 \qquad\forall n,r\in\N_0\: m,s\in\Z\]
In the classical limit this becomes the pairing $\ubp\times\cbp\to\k$
\begin{equation}
\langle X^n H^m, X^r g^s\rangle =
 \delta_{n,r} n!\, s^m\qquad \forall n,m,r\in\N_0\: s\in\Z
\label{eq:pair_class}
\end{equation}

\subsection{Differential Calculi and Quantum Tangent Spaces}

In this section we recall some facts about differential calculi
along the lines of Majid's treatment in \cite{Majid_calculi}.

Following Woronowicz \cite{Wor_calculi}, first order bicovariant differential
calculi on a quantum group $A$ (of
function algebra type) are in one-to-one correspondence to submodules
$M$ of $\ker\cou\subset A$ in the category $^A_A\cal{M}$ of (say) left
crossed modules of $A$ via left multiplication and left adjoint
coaction:
\[
a\act v = av \qquad \mathrm{Ad_L}(v)
 =v_{(1)}\antip v_{(3)}\tensor v_{(2)}
\qquad \forall a\in A, v\in A
\]
More precisely, given a crossed submodule $M$, the corresponding
calculus is given by $\Gamma=\ker\cou/M\tensor A$ with $\diff a =
\pi(\cop a - 1\tensor a)$ ($\pi$ the canonical projection).
The right action and coaction on $\Gamma$ are given by
the right multiplication and coproduct on $A$, the left action and
coaction by the tensor product ones with $\ker\cou/M$ as a left
crossed module. In all of what follows, ``differential calculus'' will
mean ``bicovariant first order differential calculus''.

Alternatively \cite{Majid_calculi}, given in addition a quantum group $H$
dually paired with $A$
(which we might think of as being of enveloping algebra type), we can
express the coaction of $A$ on
itself as an action of $H^{op}$ using the pairing:
\[
h\act v = \langle h, v_{(1)} \antip v_{(3)}\rangle v_{(2)}
\qquad \forall h\in H^{op}, v\in A
\]
Thereby we change from the category of (left) crossed $A$-modules to
the category of left modules of the quantum double $A\xprod H^{op}$.

In this picture the pairing between $A$ and $H$ descends to a pairing
between $A/\k 1$ (which we may identify with $\ker\cou\subset A$) and
$\ker\cou\subset H$. Further quotienting $A/\k 1$ by $M$ (viewed in
$A/\k 1$) leads to a pairing with the subspace $L\subset\ker\cou H$
that annihilates $M$. $L$ is called a ``quantum tangent space''
and is dual to the differential calculus $\Gamma$ generated by $M$ in
the sense that $\Gamma\cong \Lin(L,A)$ via
\begin{equation}
A/(\k 1+M)\tensor A \to \Lin(L,A)\qquad
v\tensor a \mapsto \langle \cdot, v\rangle a
\label{eq:eval}
\end{equation}
if the pairing between $A/(\k 1+M)$ and $L$ is non-degenerate.

The quantum tangent spaces are obtained directly by dualising the
(left) action of the quantum double on $A$ to a (right) action on
$H$. Explicitly, this is the adjoint action and the coregular action
\[
h \act x = h_{(1)} x \antip h_{(2)} \qquad
a \act x = \langle x_{(1)}, a \rangle x_{(2)}\qquad
 \forall h\in H, a\in A^{op},x\in A
\]
where we have converted the right action to a left action by going
from \mbox{$A\xprod H^{op}$}-modules to \mbox{$H\xprod A^{op}$}-modules.
Quantum tangent spaces are subspaces of $\ker\cou\subset H$ invariant
under the projection of this action to $\ker\cou$ via \mbox{$x\mapsto
x-\cou(x) 1$}. Alternatively, the left action of $A^{op}$ can be
converted to a left coaction of $H$ being the comultiplication (with
subsequent projection onto $H\tensor\ker\cou$).

We can use the evaluation map (\ref{eq:eval})
to define a ``braided derivation'' on elements of the quantum tangent
space via
\[\partial_x:A\to A\qquad \partial_x(a)={\diff a}(x)=\langle
x,a_{(1)}\rangle a_{(2)}\qquad\forall x\in L, a\in A\]
This obeys the braided derivation rule
\[\partial_x(a b)=(\partial_x a) b
 + a_{(2)} \partial_{a_{(1)}\act x}b\qquad\forall x\in L, a\in A\]

Given a right invariant basis $\{\eta_i\}_{i\in I}$ of $\Gamma$ with a
dual basis $\{\phi_i\}_{i\in I}$ of $L$ we have
\[{\diff a}=\sum_{i\in I} \eta_i\cdot \partial_i(a)\qquad\forall a\in A\]
where we denote $\partial_i=\partial_{\phi_i}$. (This can be easily
seen to hold by evaluation against $\phi_i\ \forall i$.)

\section{Classification on \cqbp{} and \uqbp{}}
\label{sec:q}

In this section we completely classify differential calculi on \cqbp{}
and, dually, quantum tangent spaces on \uqbp{}. We start by
classifying the relevant crossed modules and then proceed to a
detailed description of the calculi.

\begin{lem}
\label{lem:cqbp_class}
(a) Left crossed \cqbp-submodules $M\subseteq\cqbp$ by left
multiplication and left
adjoint coaction are in one-to-one correspondence to
pairs $(P,I)$
where $P\in\k(q)[g]$ is a polynomial with $P(0)=1$ and $I\subset\N$ is
finite.
$\codim M<\infty$ iff $P=1$. In particular $\codim M=\sum_{n\in I}n$
if $P=1$.

(b) The finite codimensional maximal $M$
correspond to the pairs $(1,\{n\})$ with $n$ the
codimension. The infinite codimensional maximal $M$ are characterised by
$(P,\emptyset)$ with $P$ irreducible and $P(g)\neq 1-q^{-k}g$ for any
$k\in\N_0$.

(c) Crossed submodules $M$ of finite
codimension are intersections of maximal ones.
In particular $M=\bigcap_{n\in I} M^n$, with $M^n$ corresponding to
$(1,\{n\})$.
\end{lem}
\begin{proof}
(a) Let $M\subseteq\cqbp$ be a crossed \cqbp-submodule by left
multiplication and left adjoint coaction and let
$\sum_n X^n P_n(g) \in M$, where $P_n$ are polynomials in $g,g^{-1}$
(every element of \cqbp{} can be expressed in
this form). From the formula for the coaction ((\ref{eq:adl}), see appendix)
we observe that for all $n$ and for all $t\le n$ the element
\[X^t P_n(g) \prod_{s=1}^{n-t} (1-q^{s-n}g)\]
lies in $M$.
In particular
this is true for $t=n$, meaning that elements of constant degree in $X$
lie separately in $M$. It is therefore enough to consider such
elements.

Let now $X^n P(g) \in M$.
By left multiplication $X^n P(g)$ generates any element of the form
$X^k P(g) Q(g)$, where $k\ge n$ and $Q$ is any polynomial in
$g,g^{-1}$. (Note that $Q(q^kg) X^k=X^k Q(g)$.)
We see that $M$ contains the following elements:
\[\begin{array}{ll}
\vdots & \\
X^{n+2} & P(g) \\
X^{n+1} & P(g) \\
X^n & P(g) \\
X^{n-1} & P(g) (1-q^{1-n}g) \\
X^{n-2} & P(g) (1-q^{1-n}g) (1-q^{2-n}g) \\
\vdots & \\
X & P(g) (1-q^{1-n}g) (1-q^{2-n}g) \ldots (1-q^{-1}g) \\
& P(g) (1-q^{1-n}g) (1-q^{2-n}g) \ldots (1-q^{-1}g)(1-g) 
\end{array}
\]
Moreover, if $M$ is generated by $X^n P(g)$ as a module
then these elements generate a basis for $M$ as a vector
space by left
multiplication with polynomials in $g,g^{-1}$. (Observe that the
application of the coaction to any of the elements shown does not
generate elements of new type.)

Now, let $M$ be a given crossed submodule. We pick, among the
elements in $M$ of the form $X^n P(g)$ with $P$ of minimal
length,
one
with lowest degree in $X$. Then certainly the elements listed above are
in $M$. Furthermore for any element of the form $X^k Q(g)$, $Q$ must
contain $P$ as a factor and for $k<n$, $Q$ must contain $P(g) (1-q^{1-n}g)$
as a factor. We continue by picking the smallest $n_2$, so that
$X^{n_2} P(g) (1-q^{1-n}g) \in M$. Certainly $n_2<n$. Again, for any
element of $X^l Q(g)$ in $M$ with $l<n_2$, we have that
$P(g) (1-q^{1-n}g) (1-q^{1-n_2}g)$ divides Q(g). We proceed by
induction, until we arrive at degree zero in $X$.

We obtain the following elements generating a basis for $M$ by left
multiplication with polynomials in $g,g^{-1}$ (rename $n_1=n$):

\[ \begin{array}{ll}
\vdots & \\
X^{n_1+1} & P(g) \\
X^{n_1} & P(g) \\
X^{n_1-1} & P(g) (1-q^{1-{n_1}}g) \\
\vdots & \\
X^{n_2} & P(g) (1-q^{1-{n_1}}g) \\
X^{n_2-1} & P(g) (1-q^{1-{n_1}}g) (1-q^{1-n_2})\\
\vdots & \\
X^{n_3} & P(g) (1-q^{1-{n_1}}g) (1-q^{1-{n_2}}g) \\
X^{n_3-1} & P(g) (1-q^{1-{n_1}}g) (1-q^{1-{n_2}}g) (1-q^{1-n_3})\\
\vdots & \\
& P(g) (1-q^{1-{n_1}}g) (1-q^{1-n_2}g) (1-q^{1-n_3}g) \ldots (1-q^{1-n_m}g) 
\end{array}
\]
We see that the integers $n_1,\ldots,n_m$ uniquely determine the shape
of this picture. The polynomial $P(g)$ on the other hand can be
shifted (by $g$ and $g^{-1}$) or renormalised. To determine $M$
uniquely we shift and normalise $P$ in such a way that it contains no
negative powers
and has unit constant coefficient. $P$ can then be viewed as a
polynomial $\in\k(q)[g]$.

We see that the codimension of $M$ is the sum of the lengths of the
polynomials in $g$ over all degrees in $X$ in the above
picture. Finite codimension corresponds to $P=1$. In this
case the codimension is the sum
$n_1+\ldots +n_m$.

(b) We observe that polynomials of the form $1-q^{j}g$
have no common divisors for distinct $j$. Therefore,
finite codimensional crossed
submodules are maximal if and only if
there is just one integer ($m=1$). Thus, the maximal left
crossed submodule of
codimension $k$ is generated by $X^k$ and $1-q^{1-k}g$.
For an infinite codimensional crossed submodule we certainly need
$m=0$. Then, the maximality corresponds to irreducibility of
$P$.

(c) This is again due to the distinctness of factors $1-q^j g$.
\end{proof}

\begin{cor}
\label{cor:cqbp_eclass}
(a) Left crossed \cqbp-submodules $M\subseteq\ker\cou\subset\cqbp$
are in one-to-one correspondence to pairs
$(P,I)$ as in lemma \ref{lem:cqbp_class}
with the additional constraint $(1-g)$ divides $P(g)$ or $1\in I$.
$\codim M<\infty$ iff $P=1$. In particular $\codim M=(\sum_{n\in I}n)-1$
if $P=1$.

(b) The finite codimensional maximal $M$
correspond to the pairs
$(1,\{1,n\})$ with $n\ge 2$ the
codimension. The infinite codimensional maximal $M$ correspond to pairs
$(P,\{1\})$ with $P$ irreducible and $P(g)\neq 1-q^{-k}g$ for any
$k\in\N_0$.

(c) Crossed submodules $M$ of finite
codimension are intersections of maximal ones.
In particular $M=\bigcap_{n\in I} M^n$, with $M^n$ corresponding to
$(1,\{1,n\})$.
\end{cor}
\begin{proof}
First observe that $\sum_n X^n P_n(g)\in \ker\cou$ if and only if
$(1-g)$ divides $P_0(g)$. This is to say that that $\ker\cou$
is the crossed submodule corresponding to the pair $(1,\{1\})$ in
lemma \ref{lem:cqbp_class}. We obtain the classification
from the one of lemmas \ref{lem:cqbp_class} by intersecting
everything with this crossed submodule. In particular, this reduces
the codimension by one in the finite codimensional case.
\end{proof}

\begin{lem}
\label{lem:uqbp_class}
(a) Left crossed \uqbp-submodules $L\subseteq\uqbp$ via the left adjoint
action and left
regular coaction are in one-to-one correspondence to the set
$3^{\N_0}\times2^{\N}$.
Finite dimensional $L$ are in one-to-one correspondence to
finite sets $I\subset\N$ and $\dim L=\sum_{n\in I}n$.

(b) Finite dimensional irreducible $L$ correspond to $\{n\}$
with $n$ the dimension.

(c) Finite dimensional $L$ are direct sums of irreducible ones. In
particular $L=\oplus_{n\in I} L^n$ with $L^n$ corresponding to $\{n\}$.
\end{lem}
\begin{proof}
(a) The action takes the explicit form
\[g\act X^n g^k = q^{-n} X^n g^k\qquad
X\act X^n g^k = X^{n+1}g^k(1-q^{-(n+k)})\]
while the coproduct is
\[\cop(X^n g^k)=\sum_{r=0}^{n} \binomq{n}{r}
 q^{-r(n-r)} X^{n-r} g^{k+r}\tensor X^r g^k\]
which we view as a left coaction here.
Let now $L\subseteq\uqbp$ be a crossed \uqbp-submodule via this action
and coaction. For $\sum_n X^n P_n(g)\in L$ invariance under
the action by
$g$ clearly means that \mbox{$X^n P_n(g)\in L\ \forall n$}. Then from
invariance under the coaction we can conclude that
if $X^n \sum_j a_j g^j\in L$ we must have
$X^n g^j\in L\ \forall j$.
I.e.\ elements of the form $X^n g^j$ lie separately in $L$ and it is
sufficient to consider such elements. From the coaction we learn that
if $X^n g^j\in L$ we have $X^m g^j\in L\ \forall m\le n$.
The action
by $X$ leads to $X^n g^j\in L \Rightarrow X^{n+1} g^j\in
L$ except if
$n+j=0$. The classification is given by the possible choices we have
for each power in $g$. For every positive integer $j$ we can
choose wether or not to include the span of
$\{ X^n g^j|\forall n\}$ in $L$ and for
every non-positive
integer we can choose to include either the span of $\{ X^n
g^j|\forall n\}$
or just
$\{ X^n g^j|\forall n\le -j\}$ or neither. I.e.\ for positive
integers ($\N$) we have two choices while for non-positive (identified
with $\N_0$) ones we have three choices.

Clearly, the finite dimensional $L$ are those where we choose only to
include finitely many powers of $g$ and also only finitely many powers
of $X$. The latter is only possible for the non-positive powers
of $g$.
By identifying positive integers $n$ with powers $1-n$ of $g$, we
obtain a classification by finite subsets of $\N$.

(b) Irreducibility clearly corresponds to just including one power of $g$
in the finite dimensional case.

(c) The decomposition property is obvious from the discussion.
\end{proof}

\begin{cor}
\label{cor:uqbp_eclass}
(a) Left crossed \uqbp-submodules $L\subseteq\ker\cou\subset\uqbp$ via
the left adjoint
action and left regular coaction (with subsequent projection to
$\ker\cou$ via $x\mapsto x-\cou(x)1$) are in one-to-one correspondence to
the set $3^{\N}\times2^{\N_0}$.
Finite dimensional $L$ are in one-to-one correspondence to
finite sets
$I\subset\N\setminus\{1\}$ and $\dim L=\sum_{n\in I}n$.

(b) Finite dimensional irreducible $L$ correspond to $\{n\}$
with $n\ge 2$ the dimension.

(c) Finite dimensional $L$ are direct sums of irreducible ones. In
particular $L=\oplus_{n\in I} L^n$ with $L^n$ corresponding to $\{n\}$.
\end{cor}
\begin{proof}
Only a small modification of lemma \ref{lem:uqbp_class} is
necessary. Elements of
the form $P(g)$ are replaced by elements of the form
$P(g)-P(1)$. Monomials with non-vanishing degree in $X$ are unchanged.
The choices for elements of degree $0$ in $g$ are reduced to either
including the span of
$\{ X^k |\forall k>0 \}$ in the crossed submodule or not. In
particular, the crossed submodule characterised by \{1\} in lemma
\ref{lem:uqbp_class} is projected out.
\end{proof}

Differential calculi in the original sense of Woronowicz are
classified by corollary \ref{cor:cqbp_eclass} while from the quantum
tangent space
point of view the
classification is given by corollary \ref{cor:uqbp_eclass}.
In the finite dimensional case the duality is strict in the sense of a
one-to-one correspondence.
The infinite dimensional case on the other hand depends strongly on
the algebraic models we use for the function or enveloping
algebras. It is therefore not surprising that in the present purely
algebraic context the classifications are quite different in this
case. We will restrict ourselves to the finite dimensional
case in the following description of the differential calculi.

\begin{thm}
\label{thm:q_calc}
(a) Finite dimensional differential calculi $\Gamma$ on \cqbp{} and
corresponding quantum tangent spaces $L$ on \uqbp{} are
in one-to-one correspondence to
finite sets $I\subset\N\setminus\{1\}$. In particular
$\dim\Gamma=\dim L=\sum_{n\in I}n$.

(b) Coirreducible $\Gamma$ and irreducible $L$ correspond to
$\{n\}$ with $n\ge 2$ the dimension.
Such a $\Gamma$ has a
right invariant basis $\eta_0,\dots,\eta_{n-1}$ so that the relations
\begin{gather*}
\diff X=\eta_1+(q^{n-1}-1)\eta_0 X \qquad
 \diff g=(q^{n-1}-1)\eta_0 g\\
[a,\eta_0]=\diff a\quad \forall a\in\cqbp\\
[g,\eta_i]_{q^{n-1-i}}=0\quad \forall i\qquad
[X,\eta_i]_{q^{n-1-i}}=\begin{cases}
 \eta_{i+1} & \text{if}\ i<n-1 \\
 0 & \text{if}\ i=n-1
 \end{cases}
\end{gather*}
hold, where $[a,b]_p := a b - p b a$. By choosing the dual basis on
the corresponding irreducible $L$ we obtain
the braided derivations
\begin{gather*}
\partial_i\no{f}=
 \no{Q_{n-1-i,g} Q_{n-1-i,X} \frac{1}{[i]_q!} (\partial_{q,X})^i f}
 \qquad\forall i\ge 1\\
\partial_0\no{f}=
 \no{Q_{n-1,g} Q_{n-1,X} f - f}
\end{gather*}
for $f\in \k(q)[X,g,g^{-1}]$ with normal ordering
$\k(q)[X,g,g^{-1}]\to \cqbp$ given by \mbox{$g^n X^m\mapsto g^n X^m$}.

(c) Finite dimensional $\Gamma$ and $L$ decompose into direct sums of
coirreducible respectively irreducible ones.
In particular $\Gamma=\oplus_{n\in I}\Gamma^n$ and
$L=\oplus_{n\in I}L^n$ with $\Gamma^n$ and $L^n$ corresponding to $\{n\}$.
\end{thm}
\begin{proof}
(a) We observe that the classifications of lemma
\ref{lem:cqbp_class} and lemma \ref{lem:uqbp_class} or
corollary \ref{cor:cqbp_eclass} and corollary \ref{cor:uqbp_eclass}
are dual to each other in the finite (co){}dimensional case. More
precisely, for $I\subset\N$ finite the crossed submodule $M$
corresponding to $(1,I)$ in lemma \ref{lem:cqbp_class} is the
annihilator of the crossed
submodule $L$ corresponding to $I$ in lemma \ref{lem:uqbp_class} 
and vice versa.
$\cqbp/M$ and $L$ are dual spaces with the induced pairing.
For $I\subset\N\setminus\{1\}$ finite this descends to 
$M$ corresponding to $(1,I\cup\{1\})$ in corollary
\ref{cor:cqbp_eclass} and $L$ corresponding to $I$ in corollary
\ref{cor:uqbp_eclass}.
For the dimension of $\Gamma$ observe
$\dim\Gamma=\dim{\ker\cou/M}=\codim M$.

(b) Coirreducibility (having no proper quotient) of $\Gamma$
clearly corresponds to maximality of $M$. The statement then follows
from parts (b) of corollaries
\ref{cor:cqbp_eclass} and \ref{cor:uqbp_eclass}. The formulas are
obtained by choosing the basis $\eta_0,\dots,\eta_{n-1}$ of
$\ker\cou/M$ as the equivalence classes of 
\[(g-1)/(q^{n-1}-1),X,\dots,X^{n-1}\]
The dual basis of $L$ is then given by
\[g^{1-n}-1, X g^{1-n},\dots, q^{k(k-1)} \frac{1}{[k]_q!} X^k g^{1-n},
\dots,q^{(n-1)(n-2)} \frac{1}{[n-1]_q!} X^{n-1} g^{1-n}\]

(c) The statement follows from corollaries \ref{cor:cqbp_eclass} and 
\ref{cor:uqbp_eclass} parts (c) with the observation
\[\ker\cou/M=\ker\cou/{\bigcap_{n\in I}}M^n
=\oplus_{n\in I}\ker\cou/M^n\]
\end{proof}

\begin{cor}
There is precisely one differential calculus on \cqbp{} which is
natural in the sense that it
has dimension $2$.
It is coirreducible and obeys the relations
\begin{gather*}
[g,\diff X]=0\qquad [g,\diff g]_q=0\qquad
[X,\diff X]_q=0\qquad [X,\diff g]_q=(q-1)({\diff X}) g
\end{gather*}
with $[a,b]_q:=ab-qba$. In particular we have
\begin{gather*}
\diff\no{f} = {\diff g} \no{\partial_{q,g} f} + {\diff X}
\no{\partial_{q,X} f}\qquad\forall f\in \k(q)[X,g,g^{-1}]
\end{gather*}
\end{cor}
\begin{proof}
This is a special case of theorem \ref{thm:q_calc}.
The formulas follow from (b) with $n=2$.
\end{proof}

\section{Classification in the Classical Limit}
\label{sec:class}

In this section we give the complete classification of differential
calculi and quantum tangent spaces in the classical case of \cbp{}
along the lines of the previous section.
We pay particular
attention to the relation to the $q$-deformed setting.

The classical limit \cbp{} of the quantum group \cqbp{} is
simply obtained by substituting the parameter $q$ with $1$.
The
classification of left crossed submodules in part (a) of lemma
\ref{lem:cqbp_class} remains
unchanged, as one may check by going through the proof.
In particular, we get a correspondence of crossed modules in the
$q$-deformed setting with crossed modules in the
classical setting
as a map of 
pairs $(P,I)\mapsto (P,I)$
that converts polynomials $\k(q)[g]$ to polynomials $\k[g]$ (if
defined) and leaves
sets $I$ unchanged. This is one-to-one in the finite
dimensional case.
However, we did use the distinctness of powers of $q$ in part (b) and
(c) of lemma
$\ref{lem:cqbp_class}$ and have to account for changing this. The
only place where we used it, was in observing that
factors $1-q^j g $ have no common divisors for distinct $j$. This was
crucial to conclude the maximality (b) of certain finite codimensional
crossed submodules and the intersection property (c).
Now, all those factors become $1-g$.

\begin{cor}
\label{cor:cbp_class}
(a) Left crossed \cbp-submodules $M\subseteq\cbp$ by left
multiplication and left
adjoint coaction are in one-to-one correspondence to
pairs $(P,I)$
where $P\in\k[g]$ is a polynomial with $P(0)=1$ and $I\subset\N$ is
finite.
$\codim M<\infty$ iff $P=1$. In particular $\codim M=\sum_{n\in I}n$
if $P=1$.

(b) The infinite codimensional maximal $M$ are characterised by
$(P,\emptyset)$ with $P$ irreducible and $P(g)\neq 1-g$ for any
$k\in\N_0$.
\end{cor}

In the restriction to $\ker\cou\subset\cbp$ corresponding to corollary
\ref{cor:cqbp_eclass} we observe another difference to the
$q$-deformed setting.
Since the condition for a crossed submodule to lie in $\ker\cou$ is exactly
to have factors $1-g$ in the $X$-free monomials this condition may now
be satisfied more easily. If the characterising polynomial does not
contain this factor it is now sufficient to have just any non-empty
characterising integer set $I$ and it need not contain $1$. Consequently,
the map $(P,I)\mapsto (P,I)$ does not reach all crossed submodules now.

\begin{cor}
\label{cor:cbp_eclass}
(a) Left crossed \cbp-submodules $M\subseteq\ker\cou\subset\cbp$
are in one-to-one correspondence to pairs
$(P,I)$ as in corollary \ref{cor:cbp_class}
with the additional constraint $(1-g)$ divides $P(g)$ or $I$ non-empty.
$\codim M<\infty$ iff $P=1$. In particular $\codim M=(\sum_{n\in I}n)-1$
if $P=1$.

(b) The infinite codimensional maximal $M$ correspond to pairs
$(P,\{1\})$ with $P$ irreducible and $P(g)\neq 1-g$.
\end{cor}

Let us now turn to quantum tangent spaces on \ubp{}. Here, the process
to go from the $q$-deformed setting to the classical one is not quite
so straightforward.

\begin{lem}
\label{lem:ubp_class}
Proper left crossed \ubp-submodules $L\subset\ubp$ via the left
adjoint action
and left regular coaction are
in one-to-one correspondence to pairs $(l,I)$ with $l\in\N_0$ and
$I\subset\N$ finite. $\dim L<\infty$ iff $l=0$. In particular $\dim
L=\sum_{n\in I}n$ if $l=0$.
\end{lem}
\begin{proof}
The left adjoint action takes the form
\[
X\act X^n H^m = X^{n+1}(H^m-(H+1)^m) \qquad
H\act X^n H^m = n X^n H^m
\]
while the coaction is
\[
\cop(X^n H^m) = \sum_{i=1}^n \sum_{j=1}^m \binom{n}{i} \binom{m}{j}
X^i H^j\tensor X^{n-1} H^{m-j}
\]
Let $L$ be a crossed submodule invariant under the action and coaction.
The (repeated) action of $H$ separates elements by degree in $X$. It is
therefore sufficient to consider elements of the form $X^n P(H)$, where
$P$ is a polynomial.
By acting with $X$ on an element $X^n P(H)$ we obtain
$X^{n+1}(P(H)-P(H+1))$. Subsequently applying the coaction and
projecting on the left hand side of the tensor product onto $X$ (in
the basis $X^i H^j$ of \ubp)
leads to the element $X^n (P(H)-P(H+1))$. Now the degree of
$P(H)-P(H+1)$ is exactly the degree of $P(H)$ minus $1$. Thus we have
polynomials $X^n P_i(H)$ of any degree $i=\deg(P_i)\le \deg(P)$ in $L$
by induction. In particular, $X^n H^m\in L$ for all
$m\le\deg(P)$. It is thus sufficient to consider elements of
the form $X^n H^m$. Given such an element, the coaction generates all
elements of the form $X^i H^j$ with $i\le n, j\le m$.

For given $n$, the characterising datum is the maximal $m$ so
that $X^n H^m\in L$. Due to the coaction this cannot decrease
with decreasing $n$ and due to the action of $X$ this can decrease at
most by $1$ when increasing $n$ by $1$. This leads to the
classification given. For $l\in N_0$ and $I\subset\N$ finite, the
corresponding crossed submodule
is generated by
\begin{gather*}
X^{n_m-1} H^{l+m-1}, X^{n_m+n_{m-1}-1} H^{l+m-2},\dots,
X^{(\sum_i n_i)-1} H^{l}\\
\text{and}\qquad
X^{(\sum_i n_i)+k} H^{l-1}\quad \forall k\ge 0\quad\text{if}\quad l>0
\end{gather*}
as a crossed module.
\end{proof}

For the transition from the $q$-deformed (lemma
\ref{lem:uqbp_class}) to the classical case we
observe that the space spanned by $g^{s_1},\dots,g^{s_m}$ with $m$
different integers $s_i\in\Z$ maps to the space spanned by
$1, H, \dots, H^{m-1}$ in the
prescription of the classical limit (as described in section
\ref{sec:intro_limits}). I.e.\ the classical crossed submodule
characterised by an integer $l$ and a finite set $I\subset\N$ comes
from a crossed submodule characterised by this same $I$ and additionally $l$
other integers $j\in\Z$ for which $X^k g^{1-j}$ is included. In
particular, we have a one-to-one correspondence in the finite
dimensional case.

To formulate the analogue of corollary \ref{cor:uqbp_eclass} for the
classical case is essentially straightforward now. However, as for
\cbp{}, we obtain more crossed submodules than those from the $q$-deformed
setting. This is due to the degeneracy introduced by forgetting the
powers of $g$ and just retaining the number of different powers. 

\begin{cor}
\label{cor:ubp_eclass}
(a) Proper left crossed \ubp-submodules
$L\subset\ker\cou\subset\ubp$ via the
left adjoint
action and left regular coaction (with subsequent projection to
$\ker\cou$ via $x\mapsto x-\cou(x)1$) are in one-to-one correspondence to
pairs $(l,I)$ with $l\in\N_0$ and $I\subset\N$ finite where $l\neq 0$
or $I\neq\emptyset$.
$\dim L<\infty$ iff $l=0$. In particular $\dim
L=(\sum_{n\in I}n)-1$ if $l=0$.
\end{cor}

As in the $q$-deformed setting, we give a description of the finite
dimensional differential calculi where we have a strict duality to
quantum tangent spaces.

\begin{prop}
(a) Finite dimensional differential calculi $\Gamma$ on \cbp{} and
finite dimensional quantum tangent spaces $L$ on \ubp{} are
in one-to-one correspondence to non-empty finite sets $I\subset\N$.
In particular $\dim\Gamma=\dim L=(\sum_{n\in I}) n)-1$.

The $\Gamma$ with $1\in\N$ are in
one-to-one correspondence to the finite dimensional
calculi and quantum tangent spaces of the $q$-deformed setting
(theorem \ref{thm:q_calc}(a)).

(b) The differential calculus $\Gamma$ of dimension $n\ge 2$
corresponding to the
coirreducible one of \cqbp{} (theorem \ref{thm:q_calc}(b)) has a right
invariant
basis $\eta_0,\dots,\eta_{n-1}$ so that
\begin{gather*}
\diff X=\eta_1+\eta_0 X \qquad
 \diff g=\eta_0 g\\
[g, \eta_i]=0\ \forall i \qquad
[X, \eta_i]=\begin{cases}
 0 & \text{if}\ i=0\ \text{or}\ i=n-1\\
 \eta_{i+1} & \text{if}\ 0<i<n-1
 \end{cases}
\end{gather*}
hold. The braided derivations obtained from the dual basis of the
corresponding $L$ are
given by
\begin{gather*}
\partial_i f=\frac{1}{i!}
 \left(\frac{\partial}{\partial X}\right)^i f\qquad
 \forall i\ge 1\\
\partial_0 f=\left(X \frac{\partial}{X}+
 g \frac{\partial}{g}\right) f
\end{gather*}
for $f\in\cbp$.

(c) The differential calculus of dimension $n-1$ 
corresponding to the
one in (b) with $1$ removed from the characterising set is
the same as the one above, except that we set $\eta_0=0$ and
$\partial_0=0$.
\end{prop}
\begin{proof}
(a) We observe that the classifications of corollary
\ref{cor:cbp_class} and lemma \ref{lem:ubp_class} or
corollary \ref{cor:cbp_eclass} and corollary \ref{cor:ubp_eclass}
are dual to each other in the finite (co)dimensional case.
More
precisely, for $I\subset\N$ finite the crossed submodule $M$
corresponding to $(1,I)$ in corollary \ref{cor:cbp_class} is the
annihilator of the crossed
submodule $L$ corresponding to $(0,I)$ in lemma \ref{lem:ubp_class} 
and vice versa.
$\cbp/M$ and $L$ are dual spaces with the induced pairing.
For non-empty $I$ this descends to 
$M$ corresponding to $(1,I)$ in corollary
\ref{cor:cbp_eclass} and $L$ corresponding to $(0,I)$ in corollary
\ref{cor:ubp_eclass}.
For the dimension of $\Gamma$ note
$\dim\Gamma=\dim{\ker\cou/M}=\codim M$.

(b) For $I=\{1,n\}$ we choose in
$\ker\cou\subset\cbp$ the basis $\eta_0,\dots,\eta_{n-1}$ as the
equivalence classes of
$g-1,X,\dots,X^{n-1}$. The dual basis in $L$
is then $H,X,\dots,\frac{1}{k!}X^k,\dots,\frac{1}{(n-1)!}X^{n-1}$.
This leads to the
formulas given.

(c) For $I=\{n\}$ we get the same as in (b) except that $\eta_0$ and
$\partial_0$ disappear.
\end{proof}

The classical commutative calculus is the special case of (b) with
$n=2$. It is the only calculus of dimension $2$ with
$\diff g\neq 0$. Note that it is not coirreducible.

\section{The Dual Classical Limit}
\label{sec:dual}

We proceed in this section to the more interesting point of view where
we consider the classical algebras, but with their roles
interchanged. I.e.\ we view \ubp{} as the ``function algebra''
and \cbp{} as the ``enveloping algebra''. Due to the self-duality of
\uqbp{}, we can again view the differential calculi and quantum tangent
spaces as classical limits of the $q$-deformed setting investigated in
section \ref{sec:q}.

In this dual setting the bicovariance constraint for differential
calculi becomes much
weaker. In particular, the adjoint action on a classical function
algebra is trivial due to commutativity and the adjoint coaction on a
classical enveloping algebra is trivial due to cocommutativity.
In effect, the correspondence with the
$q$-deformed setting is much weaker than in the ordinary case of
section \ref{sec:class}.
There are much more differential
calculi and quantum tangent spaces than in the $q$-deformed setting.

We will not attempt to classify all of them in the following but
essentially 
contend ourselves with those objects coming from the $q$-deformed setting.

\begin{lem}
\label{lem:cbp_dual}
Left \cbp-subcomodules $\subseteq\cbp$ via the left regular coaction are
$\Z$-graded subspaces of \cbp{} with $|X^n g^m|=n+m$,
stable under formal derivation in $X$.

By choosing any ordering in \cqbp{}, left crossed submodules via left
regular action and adjoint coaction are in one-to-one correspondence
to certain subcomodules of \cbp{} by setting $q=1$. Direct sums
correspond to direct sums.

This descends to $\ker\cou\subset\cbp$ by the projection $x\mapsto
x-\cou(x) 1$.
\end{lem}
\begin{proof}
The coproduct on \cbp{} is
\[\cop(X^n g^k)=\sum_{r=0}^{n} \binom{n}{r}
 X^{n-r} g^{k+r}\tensor X^r g^k\]
which we view as a left coaction.
Projecting on the left hand side of the tensor product onto $g^l$ in a
basis $X^n g^k$, we
observe that coacting on an element
$\sum_{n,k} a_{n,k} X^n g^k$ we obtain elements
$\sum_n a_{n,l-n} X^n g^{l-n}$ for all $l$.
I.e.\ elements of the form
$\sum_n b_n X^n g^{l-n}$ lie
separately in a subcomodule and it is
sufficient to consider such elements. Writing the coaction
on such an element as
\[\sum_t \frac{1}{t!} X^t g^{l-t}\tensor \sum_n b_n
 \frac{n!}{(n-t)!} X^{n-t} g^{l-n}\]
we see that the coaction generates all formal derivatives in $X$
of this element. This gives us the classification: \cbp-subcomodules
$\subseteq\cbp$ under the left regular coaction are $\Z$-graded
subspaces with $|X^n g^m|=n+m$, stable under formal derivation in
$X$ given by $X^n
g^m \mapsto n X^{n-1} g^m$.

The correspondence with the \cqbp case follows from
the trivial observation
that the coproduct of \cbp{} is the same as that of \cqbp{} with $q=1$.

The restriction to $\ker\cou$ is straightforward.
\end{proof}

\begin{lem}
\label{lem:ubp_dual}
The process of obtaining the classical limit \ubp{} from \uqbp{} is
well defined for subspaces and sends crossed \uqbp-submodules
$\subset\uqbp$ by
regular action and adjoint coaction to \ubp-submodules $\subset\ubp$
by regular
action. This map is injective in the finite codimensional
case. Intersections and codimensions are preserved in this case.

This descends to $\ker\cou$.
\end{lem}
\begin{proof}
To obtain the classical limit of a left ideal it is enough to
apply the limiting process (as described in section
\ref{sec:intro_limits}) to the
module generators (We can forget the additional comodule
structure). On the one hand,
any element generated by left multiplication with polynomials in
$g$ corresponds to some element generated by left multiplication with a
polynomial in $H$, that is, there will be no more generators in the
classical setting. On the other hand, left multiplication by a
polynomial in $H$ comes
from left multiplication by the same polynomial in $g-1$, that is,
there will be no fewer generators.

The maximal left crossed \uqbp-submodule $\subseteq\uqbp$
by left multiplication and adjoint coaction of
codimension $n$ ($n\ge 1$) is generated as a left ideal by
$\{1-q^{1-n}g,X^n\}$ (see lemma
\ref{lem:cqbp_class}). Applying the limiting process to this
leads to the
left ideal of \ubp{} (which is not maximal for $n\neq 1$) generated by
$\{H+n-1,X^n\}$ having also codimension $n$.

More generally, the picture given for arbitrary finite codimensional left
crossed modules of \uqbp{} in terms of generators with respect to
polynomials in $g,g^{-1}$ in lemma \ref{lem:cqbp_class} carries over
by replacing factors
$1-q^{1-n}g$ with factors $H+n-1$ leading to generators with
respect to polynomials in $H$. In particular,
intersections go to intersections since the distinctness of
the factors for different $n$ is conserved.

The restriction to $\ker\cou$ is straightforward.
\end{proof}

We are now in a position to give a detailed description of the
differential calculi induced from the $q$-deformed setting by the
limiting process.

\begin{prop}
(a) Certain finite dimensional
differential calculi $\Gamma$ on \ubp{} and quantum tangent spaces $L$
on \cbp{}
are in one-to-one correspondence to finite dimensional differential
calculi on \uqbp{} and quantum
tangent spaces on \cqbp{}. Intersections correspond to intersections.

(b) In particular,
$\Gamma$ and $L$ corresponding to coirreducible differential calculi
on \uqbp{} and
irreducible quantum tangent spaces on \cqbp{} via the limiting process
are given as follows:
$\Gamma$ has a right invariant basis
$\eta_0,\dots,\eta_{n-1}$ so that
\begin{gather*}
\diff X=\eta_1 \qquad \diff H=(1-n)\eta_0 \\
[H, \eta_i]=(1-n+i)\eta_i\quad\forall i\qquad
[X, \eta_i]=\begin{cases}
 \eta_{i+1} & \text{if}\ \ i<n-1\\
 0 & \text{if}\ \ i=n-1
\end{cases}
\end{gather*}
holds. The braided derivations corresponding to the dual basis of
$L$ are given by
\begin{gather*}
\partial_i\no{f}=\no{T_{1-n+i,H}
 \frac{1}{i!}\left(\frac{\partial}{\partial X}\right)^i f}
 \qquad\forall i\ge 1\\
\partial_0\no{f}=\no{T_{1-n,H} f - f}
\end{gather*}
for $f\in\k[X,H]$
with the normal ordering $\k[X,H]\to \ubp$ via $H^n X^m\mapsto H^n X^m$.
\end{prop}
\begin{proof}
(a) The strict duality between \cbp-subcomodules $L\subseteq\ker\cou$
given by lemma \ref{lem:cbp_dual} and corollary \ref{cor:uqbp_eclass}
and \ubp-modules $\ubp/(\k 1+M)$ with $M$ given by lemma
\ref{lem:ubp_dual} and
corollary \ref{cor:cqbp_eclass} can be checked explicitly.
It is essentially due to mutual annihilation of factors $H+k$ in
\ubp{} with elements $g^k$ in \cbp{}.

(b) $L$ is generated by
$\{g^{1-n}-1,Xg^{1-n},\dots,
X^{n-1}g^{1-n}\}$ and
$M$ is generated by $\{H(H+n-1),X(H+n-1),X^n \}$.
The formulas are obtained by denoting with
$\eta_0,\dots,\eta_{n-1}$ the equivalence classes of
$H/(1-n),X,\dots,X^{n-1}$ in $\ubp/(\k 1+M)$.
The dual basis of $L$ is then
\[g^{1-n}-1,X g^{1-n},
\dots,\frac{1}{(n-1)!}X^{n-1}
g^{1-n}\]
\end{proof}

In contrast to the $q$-deformed setting and to the usual classical
setting the many freedoms in choosing a calculus leave us with many
$2$-dimensional calculi. It is not obvious which one we should
consider to be the ``natural'' one. Let us first look at the
$2$-dimensional calculus coming from the $q$-deformed
setting as described in (b). The relations become
\begin{gather*}
[\diff H, a]=\diff a\qquad [\diff X, a]=0\qquad\forall a\in\ubp\\
\diff\no{f} =\diff H \no{\fdiff_{1,H} f} 
 + \diff X \no{\frac{\partial}{\partial X} f}
\end{gather*}
for $f\in\k[X,H]$.

We might want to consider calculi which are closer to the classical
theory in the sense that derivatives are not finite differences but
usual derivatives. Let us therefore demand
\[\diff P(H)=\diff H \frac{\partial}{\partial H} P(H)\qquad
\text{and}\qquad
\diff P(X)=\diff X \frac{\partial}{\partial X} P(X)\]
for polynomials $P$ and ${\diff X}\neq 0$ and ${\diff H}\neq 0$.

\begin{prop}
\label{prop:nat_bp}
There is precisely one differential calculus of dimension $2$ meeting
these conditions. It obeys the relations
\begin{gather*}
[a,\diff H]=0\qquad [X,\diff X]=0\qquad [H,\diff X]=\diff X\\
\diff \no{f} =\diff H \no{\frac{\partial}{\partial H} f}
 +\diff X \no{\frac{\partial}{\partial X} f}
\end{gather*}
where the normal ordering $\k[X,H]\to \ubp$ is given by
$X^n H^m\mapsto X^n H^m$.
\end{prop}
\begin{proof}
Let $M$ be the left ideal corresponding to the calculus. It is easy to
see that for a primitive element $a$ the classical derivation condition
corresponds to $a^2\in M$ and $a\notin M$. In our case $X^2,H^2\in
M$. If we take the
ideal generated from these two elements we obtain an ideal of
$\ker\cou$ of codimension $3$. Now, it is sufficient without loss of
generality to add a generator of the form $\alpha H+\beta X+\gamma
XH$. $\alpha$ and $\beta$ must then be zero in order not
to generate $X$ or $H$ in $M$.
I.e.\ $M$ is generated by $H^2,
XH, X^2$. The relations stated follow.
\end{proof}

\section{Remarks on $\kappa$-Minkowski Space and Integration}
\label{sec:kappa}

There is a straightforward generalisation of \ubp.
Let us define the Lie algebra $\lalg b_{n+}$ as generated by
$x_0,\dots, x_{n-1}$ with relations
\[ [x_0,x_i]=x_i\qquad [x_i,x_j]=0\qquad\forall i,j\ge 1\]
Its enveloping algebra \ubnp{} is nothing but (rescaled) $\kappa$-Minkowski
space as introduced in \cite{MaRu}. In this section we make some
remarks about its intrinsic geometry.

We have an injective Lie algebra
homomorphism $b_{n+}\to b_+$ given by
$x_0\mapsto H$ and $x_i\mapsto X$.
This is an isomorphism for $n=2$. The injective Lie algebra
homomorphism extends to an injective homomorphism of enveloping
algebras $\ubp\to \ubnp$ in the obvious way. This gives rise
to an injective map from the set of submodules of \ubp{} to the set of
submodules of \ubnp{} by taking the pre-image. In
particular this induces an injective
map from the set of differential calculi on \ubp{} to the set of
differential calculi on \ubnp{} which are invariant under permutations
of the $x_i\ i\ge 1$.

\begin{cor}
\label{cor:nat_bnp}
There is a natural $n$-dimensional differential calculus on \ubnp{}
induced from the one considered in proposition
\ref{prop:nat_bp}.
It obeys the relations
\begin{gather*}
[a,\diff x_0]=0\quad\forall a\in \ubnp\qquad [x_i,\diff x_j]=0
 \quad [x_0,\diff x_i]=\diff x_i\qquad\forall i,j\ge 1\\
\diff \no{f} =\sum_{\mu=0}^{n-1}\diff x_{\mu}
 \no{\frac{\partial}{\partial x_{\mu}} f}
\end{gather*}
where the normal ordering is given by
\[\k[x_0,\dots,x_{n-1}]\to \ubnp\quad\text{via}\quad
x_{n-1}^{m_{n-1}}\cdots
x_0^{m_0}\mapsto x_{n-1}^{m_{n-1}}\cdots x_0^{m_0}\]
\end{cor}
\begin{proof}
The calculus is obtained from the ideal generated by
\[x_0^2,x_i x_j, x_i x_0\qquad\forall i,j\ge 1\]
being the pre-image of
$X^2,XH,X^2$ in \ubp{}.
\end{proof}

Let us try to push the analogy with the commutative case further and
take a look at the notion of integration. The natural way to encode
the condition of translation invariance from the classical context
in the quantum group context
is given by the condition
\[(\int\tensor\id)\circ\cop a=1 \int a\qquad\forall a\in A\]
which defines a right integral on a quantum group $A$
\cite{Sweedler}.
(Correspondingly, we have the notion of a left integral.)
Let us
formulate a slightly
weaker version of this equation
in the context of a Hopf algebra $H$ dually paired with
$A$. We write
\[\int (h-\cou(h))\act a = 0\qquad \forall h\in H, a\in A\]
where the action of $H$ on $A$ is the coregular action
$h\act a = a_{(1)}\langle a_{(2)}, h\rangle$
given by the pairing.

In the present context we set $A=\ubnp$ and $H=\cbnp$. We define the
latter as a generalisation of \cbp{} with commuting
generators $g,p_1,\dots,p_{n-1}$ and coproducts
\[\cop p_i=p_i\tensor 1+g\tensor p_i\qquad \cop g=g\tensor g\]
This can be identified (upon rescaling) as the momentum sector of the
full $\kappa$-Poincar\'e algebra (with $g=e^{p_0}$).
The pairing is the natural extension of (\ref{eq:pair_class}):
\[\langle x_{n-1}^{m_{n-1}}\cdots x_1^{m_1} x_0^{k},
  p_{n-1}^{r_{n-1}}\cdots p_1^{r_1} g^s\rangle
  = \delta_{m_{n-1},r_{n-1}}\cdots\delta_{m_1,r_1} m_{n-1}!\cdots m_1!
  s^k\]
The resulting coregular
action is conveniently expressed as (see also \cite{MaRu})
\[p_i\act\no{f}=\no{\frac{\partial}{\partial x_i} f}\qquad
  g\act\no{f}=\no{T_{1,x_0} f}\]
  with $f\in\k[x_0,\dots,x_{n-1}]$.
Due to cocommutativity, the notions of left and right integral
coincide. The invariance conditions for integration become
\[\int \no{\frac{\partial}{\partial x_i} f}=0\quad
\forall i\in\{1,\dots,n-1\} 
\qquad\text{and}\qquad \int \no{\fdiff_{1,x_0} f}=0\]
The condition on the left is familiar and states the invariance under
infinitesimal translations in the $x_i$. The condition on the right states the
invariance under integer translations in $x_0$. However, we should
remember that we use a certain algebraic model of \cbnp{}. We might add,
for example, a generator $p_0$
to \cbnp{}
that is dual to $x_0$ and behaves
as the ``logarithm'' of $g$, i.e.\ acts as an infinitesimal
translation in $x_0$. We then have the condition of infinitesimal
translation invariance
\[\int \no{\frac{\partial}{\partial x_{\mu}} f}=0\]
for all $\mu\in\{0,1,\dots,{n-1}\}$.

In the present purely algebraic context these conditions do not make
much sense. In fact they would force the integral to be zero on the
whole algebra. This is not surprising, since we are dealing only with
polynomial functions which would not be integrable in the classical
case either.
In contrast, if we had for example the algebra of smooth functions
in two real variables, the conditions just characterise the usual
Lesbegue integral (up to normalisation).
Let us assume $\k=\R$ and suppose that we have extended the normal
ordering vector
space isomorphism $\R[x_0,\dots,x_{n-1}]\cong \ubnp$ to a vector space
isomorphism of some sufficiently large class of functions on $\R^n$ with a
suitable completion $\hat{U}(\lalg{b_{n+}})$ in a functional
analytic framework (embedding \ubnp{} in some operator algebra on a
Hilbert space). It is then natural to define the integration on
$\hat{U}(\lalg{b_{n+}})$ by
\[\int \no{f}=\int_{\R^n} f\ dx_0\cdots dx_{n-1}\]
where the right hand side is just the usual Lesbegue integral in $n$
real variables $x_0,\dots,x_{n-1}$. This
integral is unique (up to normalisation) in
satisfying the covariance condition since, as we have seen,
these correspond
just to the usual translation invariance in the classical case via normal
ordering, for which the Lesbegue integral is the unique solution.
It is also the $q\to 1$ limit of the translation invariant integral on
\uqbp{} obtained in \cite{Majid_qreg}.

We see that the natural differential calculus in corollary
\ref{cor:nat_bnp} is
compatible with this integration in that the appearing braided
derivations are exactly the actions of the translation generators
$p_{\mu}$. However, we should stress that this calculus is not
covariant under the full $\kappa$-Poincar\'e algebra, since it was
shown in \cite{GoKoMa} that in $n=4$ there is no such
calculus of dimension $4$. Our results therefore indicate a new
intrinsic approach to $\kappa$-Minkowski space that allows a
bicovariant
differential calculus of dimension $4$ and a unique translation
invariant integral by normal ordering and Lesbegue integration.

\section*{Acknowledgements}
I would like to thank S.~Majid for proposing this project,
and for fruitful discussions during the preparation of this paper.

\appendix
\section*{Appendix: The adjoint coaction on \uqbp}

The coproduct on $X^n$ is
\begin{gather*}
\cop(X^n)=\sum_{r=0}^{n} \binomq{n}{r} g^r X^{n-r}\tensor X^r\\
(\id\tensor\cop)\cop(X^n)=
\sum_{r=0}^n \sum_{i=0}^r \binomq{n}{r} \binomq{r}{i}
g^r X^{n-r}\tensor g^i X^{r-i}\tensor X^i
\end{gather*}
From this we get
\begin{equation*}
\begin{split}
\adl(X^n)&=
\sum_{r=0}^n \sum_{s=0}^r \binomq{n}{r} \binomq{r}{s}
g^r X^{n-r}\antip X^s\tensor g^s X^{r-s} \\
&=
\sum_{r=0}^n \sum_{s=0}^r \binomq{n}{r} \binomq{r}{s}
g^r X^{n-r} (-g^{-1}X)^s\tensor g^s X^{r-s} \\
&=
\sum_{t=0}^n \sum_{s=0}^{n-t} \binomq{n}{t+s} \binomq{t+s}{s}
g^{t+s} X^{n-t-s} (-g^{-1}X)^s\tensor g^s X^t \\
&=
\sum_{t=0}^n \sum_{s=0}^{n-t} \binomq{n}{t} \binomq{n-t}{s}
g^{t+s} X^{n-t-s} (-g^{-1}X)^s\tensor g^s X^t \\
&=
\sum_{t=0}^n \binomq{n}{t} g^t X^{n-t}
\tensor X^t \sum_{s=0}^{n-t} \binomq{n-t}{s}
q^{s(s+1)/2} (-q^{-n}g)^s \\
&=
\sum_{t=0}^n \binomq{n}{t} g^t X^{n-t}
\tensor X^t \prod_{u=1}^{n-t} (1-q^{u-n}g)\\
\end{split}
\end{equation*}
where we have used
\[\sum_{i=0}^n \binomq{n}{i} q^{i(i+1)/2} x^i=\prod_{j=1}^n (1+q^j
x)\]
which can be easily checked by induction.
Using the property
\[\adl(a g^n)=\adl(a) (1 \tensor g^n)\qquad\forall n\in\Z\]
we obtain for any polynomial $P$ in $g,g^{-1}$:
\begin{equation}
\adl(X^n P(g))=
\sum_{t=0}^n \binomq{n}{t} g^t X^{n-t}
\tensor X^t P(g) \prod_{u=1}^{n-t} (1-q^{u-n}g)
\label{eq:adl}
\end{equation}

\end{document}